\input amstex
\documentstyle{amsppt}
\magnification 1200
\NoRunningHeads
\topmatter
\title Martingale approximation of non-stationary stochastic processes
\endtitle
\author Dalibor Voln\'y \endauthor
\affil
Universit\'e de Rouen \endaffil
\address{Universit\'e de Rouen, Laboratoire Rapha\"el Salem, 
F-76821 Mont-Saint-Aignan Cedex, France}
\endaddress
\email{Dalibor.Volny\@univ-rouen.fr}
\endemail
\abstract{We generalise the martingale-coboundary representation of discrete time
stochastic processes to the non-stationary case and to random variables in Orlicz
spaces. Related limit theorems (CLT, invariance principle, log log law,
probabilities of large deviations) are studied.}
\endabstract
\endtopmatter
\document

\subheading{0. Introduction}
Let $(\Omega ,\Cal A ,\mu)$ be a probability space and $(X_i)_{i\in\Bbb Z}$ 
a sequence of random variables. We shall say that $(X_i)$ is stationary (we shall always
mean strictly stationary) if there exists a measurable and measure preserving
transformation $T$ : $\Omega\to \Omega$ such that $X_{i+1} = X_i\circ T$; each stationary
sequence can be represented in this way.
We shall study approximations of the sequence $(X_i)$ of the form
$$
  X_i = Y_i + U_i - U_{i+1}                \tag1
$$
where $(Y_i)$ is a martingale difference sequence with a filtration $(\Cal F_i)$.
In the stationary case we suppose that $\Cal F_i = T^{-i}\Cal F_0$, $i\in\Bbb Z$.
For the partial sums $S_n = \sum_{i=1}^n X_i$ we get 
$S_n = \sum_{i=0}^{n-1} Y_i + U_0 - U_n$ hence if $Y_i, U_i$ are square integrable
and the norms of $U_i$ are uniformly bounded, 
then
$$
  \frac1{\sqrt n} \Big\|\sum_{i=0}^{n-1} (X_i-Y_i)\Big\|_2 \to 0 \tag{G1}
$$  
(in the stationary case this is equivalent to Gordin's approximation criterium from
\cite{Go1}, cf\. \cite{Vo}).
The decomposition (G1) enables us to reduce the study of the central limit theorem (CLT) 
for $(X_i)$ to the study of the CLT for the martingale difference sequence
$(Y_i)$. If $(Y_i)$ and $(U_i)$ are stationary sequences with finite second moments,
(1) guarantees the Donsker invariance principle and the functional law of iterated logarithm
(cf\. \cite{Ha-He}, \cite{Vo-Sa}). Even in the general non-stationary case 
we can using (1) reduce the study of the probabilities of large deviations for $(X_i)$ 
to the study of $(Y_i)$ (cf\. \cite{Les-Vo}).

Gordin (\cite{Go1}, Theorem 2) proved (1) for a stationary sequence of
$X_i\in L^{2+\delta}$ where for $p = (2+\delta)/(1+\delta)$, $\delta\geq 0$, it is
$$
  \sum_{i=0}^\infty \|E(X_{i}|\Cal F_0)\|_p<\infty,\quad
  \sum_{i=1}^\infty \|X_{-i} - E(X_{-i}|\Cal F_0)\|_p<\infty. 
$$  
In \cite{Vo} it is proved that for a stationary sequence of $X_i\in L^p$
which are $\Cal F_\infty$ measurable and $E(X_i|\Cal F_{-\infty})=0$,
$1\leq p<\infty$ (the proof was given for $p=1,2$ but is valid for the other values of $p$
as well), the convergence of 
$$
  \sum_{i=0}^\infty E(X_{i}|\Cal F_0),\quad
  \sum_{i=1}^\infty [X_{-i} - E(X_{-i}|\Cal F_0)] 
$$ 
in $L^p$ is a necessary and sufficient condition for (1).

The idea of the martingale-coboundary decomposition was in its nonstationary form used in 
\cite{Ph-St, Lemma~2.1};
in their theorem, the random variables $X_i$ are adapted to the filtration 
$(\Cal F_i)$ and (1) follows from the convergence of the series 
$\sum_{i=0}^\infty E(X_{k+i}|\Cal F_k)$; the convergence is understood as the convergence of 
$$
  \sum_{i=0}^\infty E|E(X_{k+i}|\Cal F_k)|  
$$ 
(the $Y_i$ and $U_i$ are then integrable).

If $X_i = U_i - U_{i+1}$ and the sequence of $U_i$ is stochastically bounded then the partial 
sums $S_n$ are stochastically bounded as well and in the limit theorems we get a degenerate 
case. As shown in \cite{Sc}, in the stationary case the $S_n$ are stochastically bounded
if and only if $X_i = U_i - U_{i+1}$. The result
was generalised to nonstationary sequences in \cite{Br}.

The decomposition (1) applies to various interesting and important cases. Many of them were
studied in the monography \cite{Ha-He}, applications to mixing processes
appear in \cite{Du-G}, in \cite{Le} the method is used to limit
theorems for nonhyperbolic toral automorphisms.

\smallskip
In this paper we shall give necessary and sufficient condition of (1) in the
non-stationary case. The random variables are defined in more general function
spaces. Several limit theorems (the CLT, the invariance principle, the law of iterated logarithm, the estimates of probabilities of large deviations) implied by the decomposition 
are shown. 
\medskip

\subheading{1. Main result} By $\Cal L$ we
denote a Banach space of (classes of) integrable functions where
\roster
\item"-" the norm of $X\in \Cal L$ does not depend but on the distribution
of $X$ (hence two functions equal a.s\. are considered equal);
\item"-" $\Cal L$ is a vector space closed with respect to conditional
expectation; the operations of addition and
scalar multiplication are continuous.
\endroster

\proclaim{Theorem 1}
Let $(X_i)$, $i\in\Bbb Z$,  be a sequence of random variables, 
$X_i \in \Cal L$,
and $(\Cal F_i)$ a filtration. Suppose that all $X_i$ are
$\Cal F_\infty$-measurable and $E(X_i|\Cal F_{-\infty})=0$.

\underbar{A}. It is equivalent:
\roster
\item"(i)" There exist  a sequence of martingale differences $Y_i\in \Cal L$
with the filtration $(\Cal F_i)$ and random variables $U_i\in \Cal L$ such that
for all $i\in \Bbb Z$
$$
  X_i = Y_i + U_i - U_{i+1} \tag1
$$
and for all $k\in \Bbb Z$, $i\to \infty$,
$$
  E(U_{k+i}|\Cal F_k) \to 0,\quad U_{k-i} - E(U_{k-i}|\Cal F_k) \to 0
  \tag2
$$
in $L^1$.
 
\item"(ii)" For every $k\in \Bbb Z$ the sums
$$
  V_k = \sum_{i=0}^\infty E(X_{k+i}|\Cal F_{k-1}),\quad    
  W_k = \sum_{i=1}^\infty [X_{k-i} - E(X_{k-i}|\Cal F_{k-1})] \tag3
$$
converge in $L^1$ and the limits belong to $\Cal L$.
\endroster

\underbar{B}. Moreover, if (1) and (3) take place then
$$\gather
  U_k = V_k - W_k,\\
  Y_k = X_k - (V_k-W_k) + (V_{k+1}-W_{k+1}) =
  \sum_{i\in\Bbb Z} P_kX_{k+i}\endgather
$$
where for an integrable function $X$, 
$$
  P_i X = E(X|\Cal F_{i}) - E(X|\Cal F_{i-1}), \,\,\,\,i\in \Bbb Z.
$$
\endproclaim
 
\smallskip

\comment
\proclaim{Theorem 1.1$'$}
Let $\Cal L$ be a space where the martingale limit theorem takes place. Then 
it is equivalent:
\roster
\item"(i)" There exist  a sequence of martingale differences $Y_i\in
\Cal L$
with the filtration $(\Cal F_i)$ and random variables $U_i\in \Cal L$ such that
for all $i\in \Bbb Z$
$$
  X_i = Y_i + U_i - U_{i+1} \tag1
$$
and for all $k\in \Bbb Z$, $i\to \infty$,
$$
  E(U_{k+i}|\Cal F_k) \to 0,\quad U_{k-i} - E(U_{k-i}|\Cal F_k) \to 0
  \tag2
$$
in $\Cal L$.
 
\item"(ii)" For every $k\in \Bbb Z$ the sums
$$
  V_k = \sum_{i=0}^\infty E(X_{k+i}|\Cal F_k),\quad    
  W_k = \sum_{i=1}^\infty [X_{k-i} - E(X_{k-i}|\Cal F_k)] \tag3
$$
converge in $\Cal L$.
\endroster
\endproclaim
\endcomment
Recall that a strictly stationary sequence $(X_i)$ can always be represented using an 
automorphism
of a probability space: Let $\Omega'=\Bbb R^\Bbb Z$, $\Cal A'$ be the Borel $\sigma$-field
on $\Omega'$ and $\mu'$ the law of $(X_i)$ (i.e\. a  measure on $(\Omega',\Cal A')$).
If $f$ : $\Omega'\to \Bbb R$ is the coordinate projection $f$ : $\omega\mapsto \omega_0$ and 
$T$ is the left shift on
$\Bbb R^\Bbb Z$ (i.e\. $(T\omega)_i = \omega_{i+1}$) then $\mu'$ is $T$-invariant and
the sequence $(f\circ T^i)$ has the same distribution as $(X_i)$.\newline 
A sequence $(f\circ T^i)$ is a martingale difference sequence 
if and only if there exists an invariant $\sigma$-algebra $\Cal M$ (i.e\.
$\Cal M \subset T^{-1}\Cal M$) such that $f$ is $T^{-1}\Cal M$ measurable and
$E(f\,|\,\Cal M)=0$
(the sequence $\dots \subset T^{-i}\Cal M \subset T^{-i-1}\Cal M
\subset \dots$ is a filtration). Notice that 
$$
  E(f\,|\,\Cal M)\circ T^i = E(T^if\,|\,T^{-i}\Cal M).
$$
A sequence of $m\circ T^i$ is then
a martingale difference sequence adapted to $(T^{-i-1}\Cal M)$ if and only if
$m = E(m\,|\,T^{-1}\Cal M) - E(m\,|\,\Cal M)$ and the condition (1) becomes
$$
  f = m + g - g\circ T. \tag{$1'$}
$$
As a corollary to Theorem 1 we get the following statement.

\proclaim{Corollary 2} 
Let $f\in \Cal L$, $(T^{-i}\Cal M)$ be a filtration, $\Cal M_{-\infty}
= \underset i=-\infty \to{\overset \infty \to{\cap}} T^{-i}\Cal M$,
$\Cal M_{\infty} = \underset i=-\infty \to{\overset \infty \to{\vee}} 
T^{-i}\Cal M$, $E(f\,|\,\Cal M_{-\infty}) = 0$, and 
let $f$ be $\Cal M_{\infty}$-measurable. Then the conditions (i), (ii) are equivalent:
\roster
\item"(i)" There exist functions $m, g\in \Cal L$ such that
$(m\circ T^i)$ is a martingale difference sequence adapted to $(T^{-i}\Cal M)$
and
$$
  f = m + g - g\circ T.           \tag{$1'$}
$$

\item"(ii)" The sums
$$
  \sum_{k=0}^\infty E(f\circ T^k|\Cal M),\quad  \sum_{k=0}^\infty f\circ T^{-k} 
-  E(f\circ T^{-k}|\Cal M)]                                    \tag{$3'$}
$$
converge in $L^1$ and the limits belong to $\Cal L$.
\endroster
If $\Cal L$ equals $L^2$, ($3'$) is equivalent (cf\. \cite{Vo, Theorem 2.2}) 
to
$$
  \sum_{n=1}^\infty(\|\sum_{j=n}^\infty P_0(f\circ T^j)\|_2^2 +
  \|\sum_{j=n}^\infty P_0(f\circ T^{-j})\|_2^2)<\infty
$$
where $P_i$ denote the orthogonal projection operators
$P_ih=E(h|T^{-i}\Cal M)-E(h|T^{-i+1}\Cal M)$, $i\in \Bbb Z$, $h\in L^2$.
\endproclaim

Notice that we did not use an analog of (2). In fact, from
$E(f-m|\Cal M_{-\infty}) = 0$ it follows $E(g|\Cal M_{-\infty})=
E(g\circ T|\Cal M_{-\infty}) = E(g|\Cal M_{-\infty})\circ T$ hence we can
suppose $E(g|\Cal M_{-\infty})=0$. Similarly, because $E(f-m|\Cal M_{\infty})
=f-m$, we can suppose that $E(g|\Cal M_{\infty})=g$.
\medskip

\proclaim{Corollary 3}
If in the space $\Cal L$ the martingale limit theorem takes 
place, i.e\. for $X\in \Cal L$ and a filtration $(\Cal F_i)$, 
$E(X|\Cal F_i)\to E(X|\Cal F_\infty)$ in $\Cal L$ as $i\nearrow\infty$ and 
$E(X|\Cal F_i)\to E(X|\Cal F_{-\infty})$ in $\Cal L$ as $i\searrow-\infty$.
Then in Theorem 1 we can replace the convergence in $L^1$ 
by a convergence in $\Cal L$. In particular, this can be done if $\Cal L=L^p$
with $0<p<\infty$.
\endproclaim
\smallskip

In several important spaces the martingale convergence does not hold, namely in $L^\infty$
and in some Orlicz spaces.
As it will be shown in the applications, the Theorem 1.1 can still be used there.

The fact that the (backward) martingale limit theorem does not hold
in some Orlicz spaces and in $L^\infty$ has been well known among specialists but because I 
did not succeed to find a reference and the counterexamples are easy, I'll present them here.

Let $\psi(x) = e^{|x|}-1$ be a Young function (cf\. e.g\. \cite{N} for the definition)
and $L_\psi=\{X\,|\,\exists c>0\,\,\,\text{s.t.}\,\,\,E(\psi(X/c))<\infty\}$ be the 
corresponding Orlicz space with the norm $\|X\|_\psi= \inf \{c>0\,|\,
Ee^{|X|/c}<2\}$.

Let $\Omega = [0,1]$, $\Cal A$ be the Borel $\sigma$-algebra on $[0,1]$, $\mu$
the Lebesgue measure on $(\Omega,\Cal A)$. We partition the interval into adjacent 
subintervals $B_n$, $n=1,2,\dots$; each $B_n$ is a union of disjoint adjacent intervals
$A'_{n,k}$, $A''_{n,k}$, $k=2,3,\dots$. The lenghts of both $A'_{n,k}$ and $A''_{n,k}$ 
equal 
$$
  \mu(A'_{n,k}) = \mu(A''_{n,k}) = c\frac1{2^n}\frac1{e^kk\log^2k}
$$ 
where $c$ is such that 
$$
  2c\sum_{n=1}^\infty\sum_{k=2}^\infty \frac1{2^n}\frac1{e^kk\log^2k} = 1.
$$
We define the random variable $X$ equal $k$ on $A'_{n,k}$ and equal $-k$ on $A''_{n,k}$.
The $\sigma$-fields $\Cal F_n$ are generated by the sets $A'_{m,k}$, $A''_{m,k}$, 
$k=2,3,\dots$, $m\geq n$.

Then $(\Cal F_n)$ is a decreasing filtration, the random variable $X$ is 
$\Cal F_1$-measurable, $E(e^{|X|}) <\infty$ and $X_n= E(X|\Cal F_n) = X\Bbb I_{C_n}$
where $C_n=\cup_{k=n}^\infty B_k$. We get $X_n \to 0$ a.s\. Nevertheless, for each 
$\lambda>1$ we have 
$$
  E(e^{\lambda |X_n|}) = 2c \sum_{m=n}^\infty \frac{1}{2^m} 
  \sum_{k=2}^\infty \frac{e^{k(\lambda-1)}}{k\log^2k}  = \infty
$$
hence the Orlicz norm of $X_n$ is greater or equal to 1.

If we define the random variable $X$ equal $1$ on $A'_{n,k}$ and equal $-1$ on $A''_{n,k}$
we get a sequence of uniformly bounded random variables $X_n= E(X|\Cal F_n) = X\Bbb I_{C_n}$ 
converging to 0 a.s\. but not in the $L^\infty$ norm. 
\smallskip
Remark that by \cite{N} in any Orlicz space $\Cal L$, $E(X|\Cal F_i)\to 
E(X|\Cal F_\infty)$ in $\Cal L$ as $i\nearrow\infty$ .
\medskip

\demo{Proof of Theorem 1} 
The proof will follow the idea of the proof of Theorem~2.2 in \cite{Vo}.
\smallskip
1. Let us suppose (i),
$Y_k = E(Y_k|\Cal F_{k}) -E(Y_k|\Cal F_{k-1})$.
By (1) and (2)
$$\gather
  \sum_{i=0}^\infty E(X_{k+i}|\Cal F_{k-1}) =
  \sum_{i=0}^\infty \Big[E(U_{k+i}|\Cal F_{k-1}) - E(U_{k+i+1}|\Cal 
  F_{k-1})\Big] =\\
  \lim_{n\to\infty} \Big[E(U_k|\Cal F_{k-1}) - E(U_{k+n}|\Cal F_{k-1})\Big] =
  E(U_k|\Cal F_{k-1}),\\
  \sum_{i=1}^\infty [X_{k-i} - E(X_{k-i}|\Cal F_{k-1})] =
  \sum_{i=1}^\infty \Big[ (U_{k-i} - U_{k-i+1}) - E(U_{k-i}-
  U_{k-i+1}|\Cal F_{k-1})\Big] = \\
  \lim_{n\to\infty} \Big[ -U_{k} + E(U_{k}|\Cal F_{k-1}) + U_{k-n} - 
  E(U_{k-n}|\Cal F_{k-1})\Big] = -U_{k} + E(U_{k}|\Cal F_{k-1})  
  \endgather
$$
converge in $L^1$ and belong to $\Cal L$.
\medskip

2. Let us suppose (ii).

Recall that for $X$ integrable we denote
$$
  P_i X = E(X|\Cal F_{i}) - E(X|\Cal F_{i-1}), \,\,\,\,i\in \Bbb Z,
$$
and define (cf\. (3))
$$
 g_{k,i} = \cases \,\,\,\,\,\sum_{j=k}^\infty P_i X_j \quad &i\leq k-1,\\
  -\sum_{j=1}^\infty P_i X_{k-j}  &i\geq k,
  \endcases
$$
$i, k\in \Bbb Z$. Then
$$
  \sum_{i=1}^\infty g_{k,k-i} = \lim_{n\to\infty} \sum_{i=1}^{n}
  \sum_{j=k}^\infty
  P_{k-i} X_j = \lim_{n\to\infty} \sum_{j=k}^\infty \Big[E(X_j|\Cal F_{k-1}) -
  E(X_j|\Cal F_{k-n-1})\Big].
$$
By (3), $V_k= \sum_{j=k}^\infty E(X_j|\Cal F_{k-1})$ converge in $L^1$ and
belong to $\Cal L$ 
hence $E(V_k|\Cal F_{k-n-1})$ converge in $L^1$ to zero. Therefore the series
$$
  \sum_{i=1}^\infty g_{k,k-i} = \sum_{j=k}^\infty E(X_j|\Cal F_{k-1}) = V_k\,\,\,\,
  k\in\Bbb Z, 
$$
converge in $L^1$ and belong to $\Cal L$.

By the definition
$$\gather
  \sum_{i=k}^\infty g_{k,i} = -\lim_{n\to\infty} \sum_{i=k}^{k+n}
  \sum_{j=1}^\infty P_{i} X_{k-j} = \\
  -\lim_{n\to\infty} \sum_{j=1}^\infty
  \Big[E(X_{k-j}|\Cal F_{k+n}) - E(X_{k-j}|\Cal F_{k-1})\Big] =\\
  -\lim_{n\to\infty} \sum_{j=1}^\infty
  \Big[X_{k-j} - E(X_{k-j}|\Cal F_{k-1}) - \big(X_{k-j} - E(X_{k-j}|
  \Cal F_{k+n})\big)\Big].
  \endgather
$$
By (3), $W_k = \sum_{j=1}^\infty \big[X_{k-j} - E(X_{k-j}|
\Cal F_{k-1})\big]$ converge in $L^1$ and belong to $\Cal L$. 
\comment
For any $X\in L^1$ 
which is $\Cal F_\infty$-measurable and for any $n\geq 0$ we have 
$$
  X-E(X|\Cal F_{k+n}) = [X-E(X|\Cal F_{k})] - E[X-E(X|\Cal F_{k})|
  \Cal F_{k+n}],
$$
hence by contractivity of the conditional expectation in $L^1$
$$
  \sum_{j=1}^\infty \big[X_{k-j} - E(X_{k-j}|\Cal F_{k+n})\big]
  = W_k - E(W_k|\Cal F_{k+n}) \underset n \rightarrow \infty 
  \to{\longrightarrow} 0
$$ 
in $L^1$. \newline
\endcomment
We have
$$\gather
  \sum_{j=1}^\infty [X_{k-j}-E(X_{k-j}|\Cal F_{k+n})]\\
   = \sum_{j=1}^\infty \Big[[X_{k-j}-E(X_{k-j}|\Cal F_{k})] 
   - E[X_{k-j}-E(X_{k-j}|\Cal F_{k})|\Cal F_{k+n}]\Big]\\
   = W_k - E(W_k|\Cal F_{k+n}) \underset n \rightarrow \infty 
  \to{\longrightarrow} 0.
  \endgather
$$
Therefore the sum
$$
  \sum_{i=k}^\infty g_{k,i} = -\sum_{j=1}^\infty
  \big[X_{k-j} - E(X_{k-j}|\Cal F_{k-1})\big] = -W_k
$$
converges in $L^1$ and belongs to $\Cal L$.

We define
$$\aligned
  &U_k = \sum_{i\in\Bbb Z} g_{k,i} = \sum_{j=k}^\infty E(X_j|\Cal F_{k-1}) -
  \sum_{j=1}^\infty \big[X_{k-j} - E(X_{k-j}|\Cal F_{k-1})\big] = V_k-W_k, \\
  &Y_k = \sum_{i\in\Bbb Z} P_k X_{k+i};
  \endaligned \tag{11}
$$
because $V_k,W_k\in\Cal L$,
we have $U_k\in \Cal L$.\newline
Notice that for all $k,i\in \Bbb Z$, 
$$
  P_iU_k = g_{k,i}.
$$  

Now, we can prove (1) $X_k = Y_k + U_k - U_{k+1}$. This is equivalent
to the proof that for all $i\in \Bbb Z$,
$$
  P_i(X_k) = P_i(Y_k + U_k - U_{k+1}).
$$
Let $i<k$. Then
$$\gather
  P_i(Y_k + U_k - U_{k+1}) = P_i(U_k) - P_i(U_{k+1}) = g_{k,i} - 
  g_{k+1,i} =  \\
  \sum_{j=0}^\infty P_iX_{j+k} - \sum_{j=0}^\infty P_iX_{j+k+1} = 
  P_iX_k.
  \endgather
$$
Let $i>k$. Then
$$
  P_i(Y_k + U_k - U_{k+1}) = P_i(U_k) - P_i(U_{k+1}) =
  - \sum_{j=1}^\infty P_iX_{k-j} +  \sum_{j=1}^\infty P_iX_{k+1-j}  =
   P_iX_k.
$$
Finally,
$$
  P_k(Y_k + U_k - U_{k+1}) = \sum_{j\in \Bbb Z} P_kX_{k+j} -
  \sum_{j=1}^\infty P_kX_{k-j} -  \sum_{j=k+1}^\infty P_kX_{j} =
  P_kX_k.
$$
Because $Y_k = X_k+U_{k+1}-U_k$ and $X_k, U_k \in \Cal L$, we have $Y_k\in 
\Cal L$.  
\medskip
It remains to show that for $j\to \infty$,
$$
  E(U_{k+j}|\Cal F_k) \to 0,\quad U_{k-j} - E(U_{k-j}|\Cal F_k) \to 0
  \tag2
$$
in $L^1$.

Recall that by (11),
$$\gather
  U_{k+j} = V_{k+j} - W_{k+j}\\ 
  =\sum_{i=k+j}^\infty E(X_i|\Cal F_{k+j-1}) -
  \sum_{i=1}^\infty \big[X_{k+j-i} - E(X_{k+j-i}|\Cal F_{k+j-1})\big],\,\,\,\,
  k,j\in \Bbb Z.
  \endgather
$$

We have 
$$
  \sum_{i=k+j}^\infty E(X_{i}|\Cal F_k) \to 0,\,\,\,\,j\to \infty,
$$
in $L^1$ and 
$$
  E\Big(\sum_{i=1}^\infty \big[X_{k+j-i} - E(X_{k+j-i}|\Cal F_{k+j})\big]\Big|
  \Cal F_k\Big) = 0,\,\,\,\, j\geq 0
$$
hence $E(U_{k+j}|\Cal F_k) \to 0$ in $L^1$ for $j\to\infty$. 

From
$$
  \sum_{i=k-j}^\infty E(X_i|\Cal F_{k-j}) - E\Big(\sum_{i=k-j}^\infty E(X_i|
  \Cal F_{k-j})\Big|\Cal F_k\Big) = 0, \quad j\geq 0,
$$
and 
$$\multline
  \sum_{i=1}^\infty \big[X_{k-j-i} - E(X_{k-j-i}|\Cal F_{k-j})\big] -
  E\Big(\sum_{i=1}^\infty \big[X_{k-j-i} - E(X_{k-j-i}|\Cal F_{k-j})\big]
  \Big|\Cal F_k\Big)\\
  = \sum_{i=1}^\infty \big[X_{k-j-i} - E(X_{k-j-i}|\Cal F_{k})\big] \to 0,
  \,\,\,\,0\leq j\to \infty,
  \endmultline
$$
we deduce $U_{k-j} - E(U_{k-j}|\Cal F_k) \to 0$ in $L^1$.

\enddemo \qed

\subheading{2. Applications}
Using \cite{Les-Vo}, Theorem~1 gives us estimations of the 
probabilities of large deviations for $(X_i)$.

\proclaim{Proposition 4} Let $X_i, Y_i, U_i$, $i=1,2,\dots$, be
integrable random variables where $(Y_i)$ is a martingale difference sequence
and
$$
  X_i = Y_i + U_i - U_{i+1}. \tag1
$$
\roster
\item"(i)" If $Y_i$ and $U_i$ are uniformly bounded in $L^p$, $2\leq
p<\infty$, then for all $x>0$
$$
  P(|\sum_{i=1}^n X_i| > nx) = O(n^{-p/2}).
$$
\item"(ii)" If $\|Y_i\|_\infty \leq a$ and $\|U_i\|_\infty \leq b$ for all $i$
then
$$
  \mu(|\sum_{i=1}^n X_i| > xn) \leq \exp\Big(-n(x-b/n)^2\big/2a^2\Big).
$$
\item"(iii)" If there are $\lambda, c>0$ such that $Ee^{\lambda |Y_i|}$, 
$Ee^{c|U_i|}$ are uniformly bounded for all $i$, then for all $\epsilon>0$ 
there exists an $n_\epsilon$ s.t\. for all $n\geq n_\epsilon$
$$
  \mu(|\sum_{i=1}^n X_i| > xn) \leq \exp
  (-\frac12(1-\epsilon)\lambda^{2/3}x^{2/3}n^{1/3}).
$$
\endroster
\endproclaim

\demo{Proof} (i) is proved in \cite{Les-Vo, Corollary 4.4};
(ii) folows from Azuma's inequality (\cite{Az}) and (iii) follows
from \cite{Les-Vo, Theorem~3.2}.
\enddemo
\qed

As shown in \cite{Les-Vo, Theorem~3.3},  the estimate cannot be essentially improved.
\smallskip

\comment
Let $2\leq p<\infty$. If for all $i$
$$
  X_i = Y_i + U_i - U_{i+1}
$$
where $(Y_i)$ is a martingale difference sequence in $L^p$, and $\|Y_i\|_p
\leq M_1$, $\|U_i\|_{p/2} \leq M_2$, then for any $x>0$
$$
  \mu(|\sum_{i=1}^n X_i| > xn) \leq 
  \frac{1}{n^{p/2}}\Big((18pq^{1/2})^p\frac{M_1^p}{x^p} + 
  \frac{M_2^{p/2}}{x^{p/2}}\Big).
$$

If $\|Y_i\|_\infty \leq M_1$ and $\|U_i\|_\infty \leq M_2$ for all $i$
then
$$
  \mu(|\sum_{i=1}^n X_i| > xn) \leq \exp\Big(-n(x-M_2/n)^2\big/2M_1^2\Big).
$$

\endproclaim

\demo{Proof}
The proof of the first inequality follows from Theorem~3.5 in \cite{Les-Vo}
which establishes the inequality
$$
  \mu(|\sum_{i=1}^n Y_i| > xn) \leq 
  (18pq^{1/2})^p\frac{M_1^p}{x^p} \frac{1}{n^{p/2}} 
$$
and $\mu(|U_i|>xn) \leq \frac{M_2^{p/2}}{(xn)^{p/2}}$.

The second inequality follows from Azuma's inequality
$$
  \mu(|\sum_{i=1}^n Y_i| > xn) \leq e^{-nx^2/2M_1^2}.
$$  
\enddemo \qed
\endcomment

If $(X_i)$ is not stationary, the martingale-coboundary representation (1)
cannot guarantee central limit theorems. We shall thus not state limit 
theorems for 
$(X_i)$ but we shall show their equivalence to those for $(Y_i)$. Remark
that for $p\geq 2$, the (strict) stationarity of $(X_i)$ satisfying (1),
is a sufficient condition for the CLT, the law of iterated logarithm,
and the invariance principle.

\proclaim{Proposition 1.4} Let
$$
  X_i = Y_i + U_i - U_{i+1}
$$
where $(Y_i)$, $(U_i)$ are sequences of square integrable random variables, 
$$
  \liminf_{n\to\infty} \frac1{n} E \Big(\sum_{i=1}^n Y_i\Big)^2 > 0,
  \tag{4}
$$ 
$$
  \frac1{\sqrt n} U_n \to 0\quad \text{in}\,\,\,\,L^2.\tag{5}
$$
Then the central limit theorem takes place for $(X_i)$ if and only if it holds for $(Y_i)$.   

If, moreover, for all $\epsilon>0$
$$
  \frac1n \sum_{i=1}^n \int_{|U_i|>\epsilon\sqrt n} U_i^2\,d\mu \to 0
  \quad \text{as}\,\,\,\,n\to\infty \tag{6}
$$  
then the invariance principle
takes place for $(X_i)$ if and only if it holds for $(Y_i)$.

If there exists an $L^2$ random variable $Z$ such that
$$
  \mu(|U_n| > x) \leq \mu(|Z|>x/\sqrt{n \log\log n})\quad\text{for 
  all}\,\,\,\,x>0,\,\,\,\,n\geq 1 \tag{7}
$$
then the functional law of iterated logarithm takes place for $(X_i)$ if and
only if it holds for $(Y_i)$.
\endproclaim

A version of the Proposition 1.4 for the stationary and ergodic case
was proved by C.C\. Heyde in \cite{He} (cf\. also \cite{Ha-He, pp.142-143}). 
We shall follow the Heyde's proof.\newline
Remark that the decomposition (1$'$) $f = m + g - g\circ T$ guarantees the CLT even 
if $g$ is just measurable. Using a stronger version of ($3'$), in \cite{Go2},
a CLT  for $f, g$ integrable is proved. 
As shown in \cite{Vo-Sa}, if $g - g\circ T$
is square integrable but $g$ is just integrable, the functional limit
theorem and the law of iterated logarithm need not hold.
\smallskip

\demo{Proof of Proposition 1.4} We shall suppose the limit theorems for
$(Y_i)$ and derive them
for $(X_i)$; the proof of the other implication is the same.

Let us denote $\bar \sigma_n^2 = E \Big(\sum_{i=1}^n X_i\Big)^2$, 
$\sigma_n^2 = E \Big(\sum_{i=1}^n Y_i\Big)^2$.
>From (5) it follows that
$(1/\sqrt n) \sum_{i=1}^n (U_i-U_{i+1}) \to 0$ in the $L^2$ norm hence by (4)
$$
  \liminf_{n\to\infty}\sigma_n/\sqrt n >0\,\,\,\,\text{and}\,\,\,\,\sigma_n/
  \bar \sigma_n \to 1. \tag{8}
$$
\smallskip

For any $1\leq k\leq n$ we have
$$
  \Big| \frac1{\bar \sigma_n} \sum_{i=1}^k X_i -
  \frac1{\sigma_n} \sum_{i=1}^k Y_i \Big| \leq \frac1{\bar \sigma_n}
  |U_1-U_{k+1}| + \frac1{\sigma_n}\Big(\frac{\sigma_n}{\bar \sigma_n} - 1\Big)
  \Big|\sum_{i=1}^k Y_i \Big|. \tag{9}
$$  
Taking $k=n$ we deduce the equivalence for the CLT.  
\smallskip

Let us suppose that the invariance principle holds for $(Y_i)$. In order 
to prove it for $(X_i)$ it suffices to show (cf\. \cite{B})
$$
  \max_{1\leq k\leq n} \Big| \frac1{\sigma_n} \sum_{i=1}^k Y_i -
  \frac1{\bar \sigma_n} \sum_{i=1}^k X_i \Big| \to 0 \quad
  \text{in probability}.
$$
 
The inequality (9) remains valid for the maxima over $1\leq k\leq n$. By (6)
\comment
We have
$$\multline
  \max_{1\leq k\leq n} \Big| \frac1{\bar \sigma_n} \sum_{i=1}^k X_i -
  \frac1{\sigma_n} \sum_{i=1}^k Y_i \Big| \leq \\
  \max_{1\leq k\leq n} \frac1{\sigma_n} |U_1-U_{k+1}| +
  \max_{1\leq k\leq n} \frac1{\sigma_n} |U_1-U_{k+1}|
  \Big(\frac{\sigma_n}{\bar \sigma_n} - 1\Big) +\\
  \frac1{\sigma_n} \max_{1\leq k\leq n} \Big|\sum_{i=1}^k Y_i \Big| 
  \Big(\frac{\sigma_n}{\bar \sigma_n} - 1\Big). 
    \endmultline
$$
Because 
\endcomment
$$
  \mu(\max_{1\leq i\leq n} |U_i|/\sqrt n > \epsilon) \leq \sum_{i=1}^n
  \mu(U_i^2 >n\epsilon^2) \leq \frac1{n\epsilon^2}\sum_{i=1}^n
  \int_{|U_i|>\epsilon\sqrt n} U_i^2\,d\mu \to 0,
$$
from this and from (8) it follows that
the first term on the right side of the inequality (9) converges in probability to zero. 
By the invariance principle for $(Y_i)$
the distributions of $\frac1{\sigma_n} \max_{1\leq k\leq n} \Big|\sum_{i=1}^k 
Y_i \Big|$ converge weakly to the distributions of the maxima of the Brownian motion
process; from this and from (8) it follows that the second term 
on the right side of the inequality (9) converges to zero as well. 
\smallskip

Now, let us suppose (7). 
\comment
We denote
$$
  \phi(x) = (x \log\log x)^{1/2}\,\,\,\,(x>e).
$$
\endcomment
In order to prove the equivalence of the law of iterated logarithm for $(X_i)$
and $(Y_i)$ it suffices to show that
$$
  \Big| \max_{1\leq k\leq n} \frac1{\sqrt{n \log\log n}} \sum_{i=1}^k (Y_i - X_i) \Big| 
  \to 0\,\,\,\, a.s.,
$$
i.e\. that
$$
  \max_{1\leq k\leq n} \frac1{\sqrt{n \log\log n}} |U_k|
  \to 0\,\,\,\, a.s. \tag{10}
$$
By (7) we for every $\epsilon>0$ have
$$
  \sum_{n=1}^\infty \mu(U_n^2> \epsilon n \log\log n) \leq 
  \sum_{n=1}^\infty \mu(|Z|>\epsilon n) < \infty;
$$  
from this and from the Borel-Cantelli lemma it follows (10).
\enddemo \qed

\end